\documentclass[12pt,fleqn,leqno]{article}

\author{D.D.\ Poro\c sniuc\\}
\date{}
\title{A K\"ahler Einstein structure on the nonzero cotangent
bundle of a space form}

\begin{document}

\maketitle
\begin{abstract}
We obtain a K\"ahler Einstein structure on the nonzero cotangent
bundle of a Riemannian manifold of positive constant sectional
curvature. The obtained K\"ahler Einstein structure cannot have
constant holomorphic sectional curvature and is not locally
symmetric.

MSC 2000: 53C07, 53C15, 53C55.

Keywords and phrases: cotangent bundle, K\"ahler Einstein metric.

\end{abstract}

\vskip5mm {\large \bf Introduction} \vskip5mm

The differential geometry of the cotangent bundle $T^*M$ of a
Riemannian manifold $(M,g)$ is almost similar to that of the
tangent bundle $TM$. However, there are some differences, because
the lifts (vertical, complete, horizontal etc.) to $T^*M$ cannot
be defined just like in the case of $TM$.

In \cite{OprPor} V. Oproiu and the present author have obtained a
natural K\"ahler Einstein structure $(G,J)$ of diagonal type
induced on $T^*M$ from the Riemannian metric $g$. The obtained
K\"ahler structure on $T^*M$ depends on one essential parameter
$u$, which is a smooth function depending on the energy density
$t$ on $T^*M$. If the K\"ahler structure is Einstein, they get a
second order differential equation fulfilled by the parameter $u$.
In the case of the general solution, they have obtained that
$(T^*M,G,J)$ has constant holomorphic sectional curvature.

 In this paper we study the singular case where the parameter
$u=A \sqrt{t}$, $A \in {\bf R}$. The considered natural Riemannian
metric $G$ of diagonal type on the nonzero cotangent bundle
$T^*_0M$ is defined by using one parameter $v$ which is a smooth
function depending on the energy density $t$. The vertical
distribution $VT^*_0M$ and the horizontal distribution $HT^*_0M$
are orthogonal to each other with respect to the metric $G$.

 Next, the natural almost complex structure $J$ on
$T^*_0M$, that interchange the vertical and horizontal
distributions, depends of one essential parameter $v$.

 After that, we obtain that $G$ is Hermitian with respect to $J$ and it
 follows that the fundamental $2$-form $\phi$, associated to the
 almost Hermitian structure $(G,J)$ is the fundamental form defining the
 usual symplectic structure on $T^*_0M$, hence it is closed.

  From the integrability condition of $J$, it
follows that the base manifold $(M,g)$ must have positive constant
sectional curvature c and $A=\sqrt{2c}$. Moreover, if ~$ v >
-\sqrt{\frac{c}{2t}}$ ~then we get a K\"ahler structure on $
T^*_0M$ and this structure depends on one essential parameter v.

 In the case where the considered K\"ahler structure is Einstein we
 obtain a second order differential equation fulfilled by the
 parameter v and we have been able to find the general solution of
 that equation:
$$
v=\frac{(n-2)\sqrt{c}}{n\sqrt{2}}t^{-\frac{1}{2}}+At^{-\frac{n+1}{2}}+B,~~~~~~~~~~A,B\in
{\bf R}.
$$

  The obtained K\"ahler Einstein manifold cannot have
constant holomorphic sectional curvature and is not locally
symmetric.

In \cite{OprPap}, by using a Lagrangian on the Riemannian manifold
$(M,g)$, V. Oproiu and N.Papaghiuc have obtained the singular case
where
$$
v = u^{\prime} = \sqrt{\frac{c}{2t}}.
$$
They get a K\"ahler structure $(G,J)$ on the nonzero tangent
bundle of $(M,g)$.

 The manifolds, tensor fields and geometric objects we consider
in this paper, are assumed to be differentiable of class
$C^{\infty}$ (i.e. smooth). We use the computations in local
coordinates but many results from this paper may be expressed in
an invariant form. The well known summation convention is used
throughout this paper, the range for the indices $h,i,j,k,l,r,s$
being always${\{}1,...,n{\}}$ (see \cite{OprPap}, \cite{YanoIsh}).
We shall denote by ${\Gamma}(T^*_0M)$ the module of smooth vector
fields on $T^*_0M$.

\vskip5mm {\large \bf 1. Some geometric properties of $T^*M$}
\vskip5mm

Let $(M,g)$ be a smooth $n$-dimensional Riemannian manifold and
denote its cotangent bundle by $\pi :T^*M\longrightarrow M$.
Recall that there is a structure of a $2n$-dimensional smooth
manifold on $T^*M$, induced from the structure of smooth
$n$-dimensional manifold  of $M$. From every local chart
$(U,\varphi )=(U,x^1,\dots ,x^n)$  on $M$, it is induced a local
chart $(\pi^{-1}(U),\Phi )=(\pi^{-1}(U),q^1,\dots , q^n,$
$p_1,\dots ,p_n)$, on $T^*M$, as follows. For a cotangent vector
$p\in \pi^{-1}(U)\subset T^*M$, the first $n$ local coordinates
$q^1,\dots ,q^n$ are  the local coordinates $x^1,\dots ,x^n$ of
its base point $x=\pi (p)$ in the local chart $(U,\varphi )$ (in
fact we have $q^i=\pi ^* x^i=x^i\circ \pi, \ i=1,\dots n)$. The
last $n$ local coordinates $p_1,\dots ,p_n$ of $p\in \pi^{-1}(U)$
are the vector space coordinates of $p$ with respect to the
natural basis $(dx^1_{\pi(p)},\dots , dx^n_{\pi(p)})$, defined by
the local chart $(U,\varphi )$,\ i.e. $p=p_idx^i_{\pi(p)}$.

An $M$-tensor field of type $(r,s)$ on $T^*M$ is defined by sets
of $n^{r+s}$ components (functions depending on $q^i$ and $p_i$),
with $r$ upper indices and $s$ lower indices, assigned to induced
local charts $(\pi^{-1}(U),\Phi )$ on $T^*M$, such that the local
coordinate change rule is that of the local coordinate components
of a tensor field of type $(r,s)$ on the base manifold $M$ (see
\cite{Mok} for further details in the case of the tangent bundle).
An usual tensor field of type $(r,s)$ on $M$ may be thought of as
an $M$-tensor field of type $(r,s)$ on $T^*M$. If the considered
tensor field on $M$ is covariant only, the corresponding
$M$-tensor field on $T^*M$ may be identified with the induced
(pullback by $\pi $) tensor field on $T^*M$.

Some useful $M$-tensor fields on $T^*M$ may be obtained as
follows. Let $v,w:[0,\infty ) \longrightarrow {\bf R}$ be a smooth
functions and let $\|p\|^2=g^{-1}_{\pi(p)}(p,p)$ be the square of
the norm of the cotangent vector $p\in \pi^{-1}(U)$ ($g^{-1}$ is
the tensor field of type (2,0) having the components $(g^{kl}(x))$
which are the entries of the inverse of the matrix $(g_{ij}(x))$
defined by the components of $g$ in the local chart $(U,\varphi
)$). The components $g_{ij}(\pi(p))$, $p_i$, $v(\|p\|^2)p_ip_j $
define $M$-tensor fields of types $(0,2)$, $(0,1)$, $(0,2)$ on
$T^*M$, respectively. Similarly, the components $g^{kl}(\pi(p))$,
$g^{0i}=p_hg^{hi}$, $w(\|p\|^2)g^{0k}g^{0l}$ define $M$-tensor
fields of type $(2,0)$, $(1,0)$, $(2,0)$ on $T^*M$, respectively.
Of course, all the components considered above are in the induced
local chart $(\pi^{-1}(U),\Phi)$.

 The Levi Civita connection $\dot \nabla $ of $g$ defines a direct
sum decomposition

\begin{equation}
TT^*M=VT^*M\oplus HT^*M.
\end{equation}
of the tangent bundle to $T^*M$ into vertical distributions
$VT^*M= {\rm Ker}\ \pi _*$ and the horizontal distribution
$HT^*M$.

 If $(\pi^{-1}(U),\Phi)=(\pi^{-1}(U),q^1,\dots ,q^n,p_1,\dots ,p_n)$
is a local chart on $T^*M$, induced from the local chart
$(U,\varphi )= (U,x^1,\dots ,x^n)$, the local vector fields
$\frac{\partial}{\partial p_1}, \dots , \frac{\partial}{\partial
p_n}$ on $\pi^{-1}(U)$ define a local frame for $VT^*M$ over $\pi
^{-1}(U)$ and the local vector fields $\frac{\delta}{\delta
q^1},\dots ,\frac{\delta}{\delta q^n}$ define a local frame for
$HT^*M$ over $\pi^{-1}(U)$, where
$$
\frac{\delta}{\delta q^i}=\frac{\partial}{\partial
q^i}+\Gamma^0_{ih} \frac{\partial}{\partial p_h},\ \ \ \Gamma
^0_{ih}=p_k\Gamma ^k_{ih}
 $$
and $\Gamma ^k_{ih}(\pi(p))$ are the Christoffel symbols of $g$.

The set of vector fields $(\frac{\partial}{\partial p_1},\dots
,\frac{\partial}{\partial p_n}, \frac{\delta}{\delta q^1},\dots
,\frac{\delta}{\delta q^n})$ defines a local frame on $T^*M$,
adapted to the direct sum decomposition (1).

We consider
\begin{equation}
t=\frac{1}{2}\|p\|^2=\frac{1}{2}g^{-1}_{\pi(p)}(p,p)=\frac{1}{2}g^{ik}(x)p_ip_k,
\ \ \ p\in \pi^{-1}(U)
\end{equation}
the energy density defined by $g$ in the cotangent vector $p$. We
have $t\in [0,\infty)$ for all $p\in T^*M$.

From now on we shall work in a fixed local chart $(U,\varphi)$ on
$M$ and in the induced local chart $(\pi^{-1}(U),\Phi)$ on $T^*M$.

\vskip5mm {\large \bf 2. An almost K\"ahler structure on the
$T^*_0M$} \vskip5mm

The nonzero cotangent bundle $T^*_0M$ of Riemannian manifold
$(M,g)$ is defined by the formula: $T^*M$ minus zero section.
 Consider a real valued smooth
function $v$ defined on $(0,\infty)\subset {\bf R}$ and a real
constant A. We define the following $M$-tensor field of type
$(0,2)$  on $T^*_0M$ having the components
\begin{equation}
G_{ij}(p)= A\sqrt{t} g_{ij}(\pi(p))+v(t)p_ip_j.
\end{equation}

It follows easily that the matrix $(G_{ij})$ is positive definite
if and only if
\begin{equation}
A>0,\ \ \ v>-\frac{A}{2\sqrt{t}}.
\end{equation}

 The inverse of this
matrix has the entries

\begin{equation}
H^{kl}(p)= \frac{1}{A\sqrt{t}} g^{kl}(\pi(p))+w(t)g^{0k}g^{0l}.
\end{equation}
where
\begin{equation}
w(t)=\frac{-v}{At(A+2\sqrt{t}v)}.
\end{equation}
 The components $H^{kl}$ define an $M$-tensor field of type
$(2,0)$
on $T^*_0M$ .\\

{\bf Remark.} If the matrix $(G_{ij})$ is positive definite then
its inverse $(H^{kl})$ is positive definite too.

Using the $M$-tensor fields defined by $G_{ij},\ H^{kl}$, the
following Riemannian metric may be considered on $T^*_0M$:
\begin{equation}
G=G_{ij}dq^idq^j+H^{ij}Dp_iDp_j,
\end{equation}
where $Dp_i=dp_i-\Gamma^0_{ij}dq^j$ is the absolute (covariant)
differential of $p_i$ with respect to the Levi Civita connection
$\dot\nabla$ of $g$ . Equivalently, we have
$$
G(\frac{\delta}{\delta q^i},\frac{\delta}{\delta
q^j})=G_{ij},~~~G(\frac{\partial}{\partial p_i}
,\frac{\partial}{\partial p_j})=H^{ij},~~
G(\frac{\partial}{\partial p_i},\frac{\delta}{\delta q^j})=~
G(\frac{\delta}{\delta q^j},\frac{\partial}{\partial p_i})=0.
$$

Remark that $HT^*_0M,~VT^*_0M$ are orthogonal to each other with
respect to $G$, but the Riemannian metrics induced from $G$ on
$HT^*_0M,~VT^*_0M$ are not the same, so the considered metric $G$
on $T^*_0M$ is not a metric of Sasaki type.
 Remark also that the system of 1-forms
$(dq^1,...,dq^n,Dp_1,...,Dp_n)$ defines a local frame on
$T^{*}T^*_0M$, dual to the local frame $(\frac{\delta}{\delta
q^1},\dots ,\frac{\delta}{\delta q^n},\frac{\partial}{\partial
p_1},\dots ,\frac{\partial}{\partial p_n})$ adapted to the direct
sum decomposition (1).

Next, an almost complex structure $J$ is defined on $T^*_0M$ by
the same $M$-tensor fields $G_{ij},\ H^{kl}$, expressed in the
adapted local frame by
\begin{equation}
J\frac{\delta}{\delta
q^i}=G_{ik}\frac{\partial}{\partial p_k},\ \ \
J\frac{\partial}{\partial p_i}=-H^{ik}\frac{\delta}{\delta q^k}
\end{equation}

 From the property
of the $M$-tensor field $H^{kl}$ to be defined by the inverse of
the matrix defined by the components of the $M$-tensor field
$G_{ij}$, it follows easily that $J$ is an almost complex
structure on $T^*_0M$. \\

{\bf Theorem 1}. \it $(T^*_0M,G,J)$ is an almost K\"ahler
manifold.\\

 Proof. \rm Since the matrix $(H^{kl})$ is the inverse of the
 matrix $(G_{ij})$, it follows easily that
$$
G(J\frac{\delta}{\delta q^i},J\frac{\delta}{\delta
q^j})=G(\frac{\delta}{\delta q^i},\frac{\delta}{\delta q^j}),\ \ \
G(J\frac{\partial}{\partial p_i},J\frac{\partial}{\partial
p_j})=G(\frac{\partial}{\partial p_i},\frac{\partial}{\partial
p_j}),$$ $$ G(J\frac{\partial}{\partial p_i},J\frac{\delta}{\delta
q^j})=G(\frac{\partial}{\partial p_i},\frac{\delta}{\delta
q^j})=0.
$$

Hence
$$G(JX,JY)=G(X,Y),\ \ \forall\ X,Y{\in}{\Gamma}(T^*_0M).$$
Thus  $(T^*_0M,G,J)$ is an almost Hermitian manifold.

The fundamental $2$-form associated with this almost Hermitian
structure is $\phi$, defined by
$$\phi(X,Y) = G(X,JY),\ \ \ \forall\ X,Y{\in}{\Gamma}(T^*_0M).$$
By a straightforward computation we get
$$
\phi(\frac{\delta}{\delta q^i},\frac{\delta}{\delta q^j})=0,\ \ \
\phi(\frac{\partial}{\partial p_i},\frac{\partial}{\partial
p_j})=0,\ \ \ \phi(\frac{\partial}{\partial
p_i},\frac{\delta}{\delta q^j})= \delta^i_j.
$$
Hence
\begin{equation}
\phi =Dp_i\wedge dq^i= dp_i\wedge dq^i,
\end{equation}
due to the symmetry of $\Gamma^0_{ij}=p_h\Gamma^h_{ij}$. It
follows that $\phi$ does coincide with the fundamental $2$-form
defining the usual symplectic structure on $T^*_0M$. Of course, we
have $d\phi =0$, i.e. $\phi$ is closed. Therefore $(T^*_0M,G,J)$
is an almost K\"ahler manifold.

\vskip5mm {\large \bf 3. A K\"ahler structure on $T^*_0M$}
\vskip5mm

We shall study the integrability of the almost complex structure
defined by $J$ on $T^*_0M$. To do this we need the following well
known formulas for the brackets of the vector fields
$\frac{\partial}{\partial p_i},\frac{\delta}{\delta q^i},~
i=1,...,n$:
\begin{equation}
[\frac{\partial}{\partial p_i},\frac{\partial}{\partial
p_j}]=0;~~~[\frac{\partial}{\partial p_i},\frac{\delta}{\delta
q^j}]=\Gamma^i_{jk}\frac{\partial}{\partial p_k};~~~
[\frac{\delta}{\delta q^i},\frac{\delta}{\delta q^j}]
=R^0_{kij}\frac{\partial}{\partial p_k},
\end{equation}
where $R^h_{kij}(\pi(p))$ are the local coordinate components of
the curvature tensor field of $\dot \nabla$ on $M$ and
$R^0_{kij}(p)=p_hR^h_{kij}$ . Of course, the components
 $R^0_{kij}$, $R^h_{kij}$ define M-tensor fields of types
 (0,3), (1,3) on $T^*_0M$, respectively.\\

\newpage
{\bf Theorem 2. } {\it The Nijenhuis tensor field of the almost
complex structure $J$ on $T^*_0M$ is given by}
\begin{equation}
\left\{
\begin{array}{l}
N(\frac{\delta}{\delta q^i},\frac{\delta}{\delta
q^j})=\{\frac{A^2}{2}(\delta^h_ig_{jk}-
\delta^h_jg_{ik})-R^h_{kij}\}p_h\frac{\partial}{\partial p_k},
\\ \mbox{ } \\
N(\frac{\delta}{\delta q^i},\frac{\partial}{\partial
p_j})=H^{kl}H^{jr}\{\frac{A^2}{2}(\delta^h_ig_{rl}-
\delta^h_rg_{il})-R^h_{lir}\}p_h\frac{\delta}{\delta q^k},
\\ \mbox{ } \\
N(\frac{\partial}{\partial p_i},\frac{\partial}{\partial
p_j})=H^{ir}H^{jl}\{\frac{A^2}{2}(\delta^h_lg_{rk}-
\delta^h_rg_{lk})-R^h_{klr}\}p_h\frac{\partial}{\partial p_k}.
\end{array}
\right.
\end{equation}

\vskip2mm

{\it Proof. } Recall that the Nijenhuis tensor field $N$ defined
by $J$ is given by
$$
N(X,Y)=[JX,JY]-J[JX,Y]-J[X,JY]-[X,Y],\ \ \forall\ \ X,Y \in \Gamma
(T^*_0M).
$$
Then, we have $\frac{\delta}{\delta q^k}t =0,\
\frac{\partial}{\partial p_k}t = g^{0k}$ and $\dot
\nabla_iG_{jk}=0,\dot \nabla_iH^{jk}= 0$, where
$$
\dot \nabla_iG_{jk}= \frac{\delta}{\delta
q^i}G_{jk}-\Gamma^l_{ij}G_{lk}-\Gamma^l_{ik}G_{lj}
$$
$$
\dot \nabla_iH^{jk}= \frac{\delta}{\delta
q^i}H^{jk}+\Gamma^j_{il}H^{lk}+\Gamma^k_{il}H^{lj}
$$

The above expressions for the components of $N$ can be obtained by
a quite long, straightforward  computation.\\

{\bf Theorem 3.} {\it The almost complex structure $J$ on $T^*_0M$
is integrable if and only if the base manifold $M$ has positive
constant sectional curvature $c$ and}
\begin{equation}
~~~~~~~~~~~~~~~~~~~~~~~~~~~~~A=\sqrt{2c}.
\end{equation}

{\it Proof.} From the condition $N=0$, one obtains
$$
\{\frac{A^2}{2}(\delta^h_ig_{jk}-
\delta^h_jg_{ik})-R^h_{kij}\}p_h=0.
$$
Differentiating with respect to $p_l$, it follows that the
curvature tensor field of $\dot \nabla$ has the expression

$$
R^l_{kij}=\frac{A^2}{2}({\delta}^l_ig_{jk}-{\delta}^l_jg_{ik}).
$$
Thus $(M,g)$ has positive constant sectional curvature
$c=\frac{A^2}{2}$. It follows that $A=\sqrt{2c}> 0$.

 Conversely, if $(M,g)$ has positive constant
sectional curvature $c$ and $A$ is given by (11), one obtains
 in a straightforward way that $N = 0$.\\

\bf Remark. \rm  If $A=\sqrt{2c}$ then the condition (4) is
equivalent with
\begin{equation}
v>-\sqrt{\frac{c}{2t}}.
\end{equation}

Then $(T^*_0M,G,J)$ is a K\"ahler manifold.

The components of the K\"ahler metric $G$ on $T^*_0M$ are

\begin{equation}
\left\{
\begin{array}{l}
G_{ij}(p)= \sqrt{2ct} g_{ij}+vp_ip_j,
\\ \mbox{ } \\
 H^{kl}(p)=\frac{1}{\sqrt{2ct}}
g^{kl}-\frac{v}{2t(c+\sqrt{2ct}v)}g^{0k}g^{0l}.
\end{array}
\right.
\end{equation}

\vskip5mm {\large \bf 4. The Levi Civita connection of the
 metric G and its curvature tensor field}
 \vskip5mm
The Levi Civita connection $\nabla$ of the Riemannian manifold
$(T^*_0M,G)$ is determined by the conditions
$$
\nabla G=0,~~~~~  T =0,
$$
where $T$ is its torsion tensor field. The explicit expression of
this connection is obtained from the formula
$$
2G({\nabla}_XY,Z)=X(G(Y,Z))+Y(G(X,Z))-Z(G(X,Y))+
$$
$$
+G([X,Y],Z)-G([X,Z],Y)-G([Y,Z],X); ~~~~~~ \forall\
X,Y,Z~{\in}~{\Gamma}(T^*_0M).
$$

The final result can be stated as follows. \\

\bf Theorem 4. {\it The Levi Civita connection ${\nabla}$ of ~$G$
has the following expression in the local adapted frame
$(\frac{\delta}{\delta q^1},\dots ,\frac{\delta}{\delta
q^n},\frac{\partial}{\partial p_1},\dots ,\frac{\partial}{\partial
p_n}):$
$$
\left\{
\begin{array}{l}
 \nabla_\frac{\partial}{\partial
p_i}\frac{\partial}{\partial p_j}
=Q^{ij}_h\frac{\partial}{\partial p_h},\ \ \ \ \ \
\nabla_\frac{\delta}{\delta q^i}\frac{\partial}{\partial
p_j}=-\Gamma^j_{ih}\frac{\partial}{\partial
p_h}+P^{hj}_i\frac{\delta}{\delta q^h},
\\ \mbox{ } \\
\nabla_\frac{\partial}{\partial p_i}\frac{\delta}{\delta
q^j}=P^{hi}_j\frac{\delta}{\delta q^h},\ \ \ \ \ \
\nabla_\frac{\delta}{\delta q^i}\frac{\delta}{\delta
q^j}=\Gamma^h_{ij}\frac{\delta}{\delta
q^h}+S_{hij}\frac{\partial}{\partial p_h},
\end{array}
\right.
$$
where $Q^{ij}_h, P^{hi}_j, S_{hij}$ are $M$-tensor fields on
$T^*_0M$, defined by
$$
\left\{
\begin{array}{l}
 Q^{ij}_h = \frac{1}{2}G_{hk}(\frac{\partial}{\partial
p_i}H^{jk}+ \frac{\partial}{\partial p_j}H^{ik}
-\frac{\partial}{\partial p_k}H^{ij}),
\\ \mbox{ } \\
 P^{hi}_j=\frac{1}{2}H^{hk}(\frac{\partial}{\partial
p_i}G_{jk}-H^{il}R^0_{ljk}),
\\ \mbox{ } \\
S_{hij}=-\frac{1}{2}G_{hk}\frac{\partial}{\partial
p_k}G_{ij}+\frac{1}{2}R^0_{hij}.
\end{array}
\right.
$$

\rm In the case of the K\"ahler structure on $T^*_0M$, using by
the relations (14), we obtain

$$
\left\{
\begin{array}{l}
Q^{ij}_h
=-\frac{1}{4t}(\delta^i_hg^{0j}+\delta^j_hg^{0i})+\frac{c-\sqrt{2ct}~
v}{4ct}g^{ij}p_h+
\frac{v^2-\sqrt{2ct}~v^{\prime}}{4t(c+\sqrt{2ct}~v)}g^{0i}g^{0j}p_h
,
\\ \mbox{ } \\
P^{hi}_j=-Q^{ih}_j,
\\ \mbox{ } \\
S_{hij}=-\frac{c+\sqrt{2ct}~v}{2}(g_{ij}p_h+g_{hi}p_j)+\frac{c-\sqrt{2ct}~v}{2}g_{hj}p_i-
\frac{2v^2+\sqrt{2ct}~v^{\prime}+2tvv^{\prime}}{2}p_hp_ip_j.
\end{array}
\right.
$$

The curvature tensor field $K$ of the connection $\nabla $ is
obtained from the well known formula
$$
K(X,Y)Z=\nabla_X\nabla_YZ-\nabla_Y\nabla_XZ-\nabla_{[X,Y]}Z,\ \ \
\ \forall\ X,Y,Z\in \Gamma (T^*_0M).
$$

The components of curvature tensor field $K$ with respect to the
adapted local frame $(\frac{\delta}{\delta q^1},\dots
,\frac{\delta}{\delta q^n},\frac{\partial}{\partial p_1},\dots
,\frac{\partial}{\partial p_n})$ are obtained easily:
$$
\left\{
\begin{array}{l}
K(\frac{\delta}{\delta q^i},\frac{\delta}{\delta
q^j})\frac{\delta}{\delta q^k}=QQQ^h_{ijk}\frac{\delta}{\delta
q^h}=(R^h_{kij}-P^{hl}_kR^0_{lij}+P^{hl}_iS_{ljk}-P^{hl}_jS_{lik})\frac{\delta}{\delta
q^h},
\\ \mbox{ } \\
K(\frac{\delta}{\delta q^i},\frac{\delta}{\delta
q^j})\frac{\partial}{\partial
p_k}=QQP^k_{ijh}\frac{\partial}{\partial
p_h}=-QQQ^k_{ijh}\frac{\partial}{\partial p_h},
\\ \mbox{ } \\
K(\frac{\partial}{\partial p_i},\frac{\partial}{\partial
p_j})\frac{\delta}{\delta q^k}=PPQ^{ijh}_k\frac{\delta}{\delta
q^h}=(\frac{\partial}{\partial
p_i}P^{hj}_k-\frac{\partial}{\partial
p_j}P^{hi}_k+P^{hi}_lP^{lj}_k-P^{hj}_lP^{li}_k)\frac{\delta}{\delta
q^h},
\\ \mbox{ } \\
K(\frac{\partial}{\partial p_i},\frac{\partial}{\partial
p_j})\frac{\partial}{\partial p_k}
=PPP^{ijk}_h\frac{\partial}{\partial
p_h}=-PPQ^{ijk}_h\frac{\partial}{\partial p_h},
\\ \mbox{ } \\
K(\frac{\partial}{\partial p_i},\frac{\delta}{\delta
q^j})\frac{\delta}{\delta q^k}=PQQ^i_{jkh}\frac{\partial}{\partial
p_h}=(\frac{\partial}{\partial
p_i}S_{hjk}+Q^{il}_hS_{ljk}-P^{li}_kS_{hjl})\frac{\partial}{\partial
p_h},
\\ \mbox{ } \\
 K(\frac{\partial}{\partial p_i},\frac{\delta}{\delta
q^j})\frac{\partial}{\partial p_k}=PQP^{ikh}_j\frac{\delta}{\delta
q^h}=(\frac{\partial}{\partial
p_i}P^{hk}_j+P^{hi}_lP^{lk}_j-P^{hl}_jQ^{ik}_l)\frac{\partial}{\partial
p_k},
\end{array}
\right.
$$
where $QQQ^h_{ijk}$, $QQP^k_{ijh}$, $PPQ^{ijh}_k$, $PPP^{ijk}_h$,
$PQQ^i_{jkh}$, $PQP^{ikh}_j$ are M-tensor fields on $T^*_0M$. The
explicit expressions of these components are obtained after some
quite long and hard computations, made by using the package RICCI.
\\

The Ricci tensor field Ric of $\nabla$ is defined by the formula:
$$
Ric(Y,Z)=trace(X\longrightarrow K(X,Y)Z),\ \ \ \forall\  X,Y,Z\in
\Gamma (T^*_0M).
$$

It follows
$$
\left\{
\begin{array}{l}
Ric(\frac{\delta}{\delta q^j},\frac{\delta}{\delta
q^k})=RicQQ_{jk}=QQQ_{hjk}^h+PQQ_{jkh}^h,
\\ \mbox{ } \\
 Ric(\frac{\partial}{\partial
p_j},\frac{\partial}{\partial
p_k})=RicPP^{jk}=PPP^{hjk}_h-PQP^{jkh}_h,
\\ \mbox{ } \\
Ric(\frac{\partial}{\partial p_i},\frac{\delta}{\delta
q^j})=Ric(\frac{\delta}{\delta q^j},\frac{\partial}{\partial
p_i})=0.
\end{array}
\right.
$$
Doing the necessary computations, we obtain the final expressions
of the components of the Ricci tensor field of $\nabla$

\begin{equation}
\left\{
\begin{array}{l}
RicQQ_{jk}=\frac{a}{2}g_{jk}+\frac{\alpha}{4t}p_jp_k,
\\ \mbox{ } \\
RicPP^{jk}=\frac{a}{4ct}g^{jk}+\frac{\beta}{8\sqrt{c}t^2(\sqrt{c}+\sqrt{2t}~v)}g^{0j}g^{0k},
\end{array}
\right.
\end{equation}
where the coefficients $a,\alpha,\beta$ are given by
\begin{equation}
\left\{
\begin{array}{l}
a=n\sqrt{c}(\sqrt{c}-\sqrt{2}~t^\frac{1}{2}v)-\sqrt{2c}(\sqrt{2c}+t^\frac{1}{2}v+
2t^{\frac{3}{2}}v^{\prime}),
\\ \mbox{ } \\
\alpha=-n(c+2tv^2+2\sqrt{2c}~t^\frac{3}{2}v^{\prime}+4t^2vv^{\prime})+
\\ \mbox{ } \\
~~~~~+2(c-tv^2-3\sqrt{2c}~t^\frac{3}{2}v^{\prime}-8t^2vv^{\prime}-
2\sqrt{2c}~t^\frac{5}{2}v^{\prime\prime}-4t^3vv^{\prime\prime}),
\\ \mbox{ } \\
\beta=-n(c+\sqrt{2c}~t^\frac{1}{2}v-2tv^2+2\sqrt{2c}~t^\frac{3}{2}v^{\prime})+
\\ \mbox{ } \\
~~~~~+2(c+\sqrt{2c}~t^\frac{1}{2}v+tv^2-3\sqrt{2c}~t^\frac{3}{2}v^{\prime}+
2t^2vv^{\prime}-2\sqrt{2c}~t^\frac{5}{2}v^{\prime\prime}).
\end{array}
\right.
\end{equation}

\newpage
\vskip5mm {\large \bf 5. A K\"ahler Einstein structure on
$T^*_0M$} \vskip5mm

In order to find out the conditions under which the K\"ahler
structure $(T^*_0M,G,J)$ is Einstein, we consider the differences
$$
\left\{
\begin{array}{l}
DiffQQ_{jk}=RicQQ_{jk}-\frac{a}{2\sqrt{2ct}}G_{jk},
\\ \mbox{ } \\
DiffPP^{jk}=RicPP^{jk}-\frac{a}{2\sqrt{2ct}}H^{jk}.
\end{array}
\right.
$$
Using by the relations (14), (15), (16) we obtain
$$
\left\{
\begin{array}{l}
DiffQQ_{jk}=\frac{\sqrt{c}+\sqrt{2t}}{4t}~ \gamma ~ p_jp_k,
\\ \mbox{ } \\
DiffPP^{jk}=\frac{1}{8t^2(\sqrt{c}+\sqrt{2t}~v)}~\gamma~g^{0j}g^{0k},
\end{array}
\right.
$$
where the factor $\gamma$ is given by

$$
\gamma=(2-n)\sqrt{c}-2(n+3)\sqrt{2}~t^\frac{3}{2}v^{\prime}-4\sqrt{2}~t^\frac{5}{2}v^{\prime\prime}.
$$

Our purpose is to solve the system
$$
DiffQQ_{jk}=DiffPP^{jk}=0,
$$
which is equivalent with the equation
$$
\gamma=0.
$$

Finally, we have to solve the following second order differential
equation of Euler type
$$
t^2v^{\prime\prime}+\frac{n+3}{2}tv^{\prime}=\frac{(2-n)\sqrt{c}}{4\sqrt{2}}t^{-\frac{1}{2}}.
$$

The general solution of this equation is
\begin{equation}
~~~~~~~~v=\frac{(n-2)\sqrt{c}}{n\sqrt{2}}t^{-\frac{1}{2}}+At^{-\frac{n+1}{2}}+B,~~~~~~~~~~A,B\in
{\bf R}.
\end{equation}

The function v must fulfill the condition $v>-\sqrt{\frac{c}{2t}}$
which is equivalent with
\begin{equation}
A \geq 0,~~~ B \geq 0;~~~n \geq 2.
\end{equation}

From the relations (15), (16), (17) we obtain
$$
Ric=-\frac{B(n+1)}{2}G.
$$

Now we may state our main result.
 \vskip5mm

\bf Theorem 4. {\it Assume that $(M,g)$ has positive constant
sectional curvature c, the almost complex structure $J$ is given
by (8), the components of the metric G are given by (14), the
function v has the expression (17) and the conditions (18) are
fulfilled.

Then $(T^*_0M,G,J)$ is a K\"ahler Einstein manifold.} \vskip5mm

{\bf Remark.} \rm After some long and hard computations we have
obtained that the K\"ahler Einstein manifold $(T^*_0M,G,J)$ cannot
have constant holomorphic sectional curvature and is not locally
symmetric.

\vskip 1.5cm

\begin{minipage}{2.5in}
\begin{flushleft}
D.D.Poro\c sniuc\\
Department of Mathematics National College "M. Eminescu" \\
Str. Octav Onicescu 52 RO-710096 Boto\c sani, Romania.\\
e-mail: dporosniuc@yahoo.com \\
~~~~~~~~~~danielporosniuc@lme.ro
\end{flushleft}
\end{minipage}

\end{document}